\documentclass[leqno,12pt]{article} 
\usepackage{amsmath, amssymb}

\usepackage{amsthm} 
\def\A{{{\mathbb A}}}

\def\Z{{{\mathbb Z}}}
\def\Q{{{\mathbb Q}}}

\def\C{{{\mathbb C}}}
\def\P{{{\mathbb P}}}
\def\G{{{\mathbb G}}}

\def\O{{{\mathcal O}}}
\def\1C{{{\mathcal C}}}

\renewcommand\div{{{\text{div}}}}

\theoremstyle{plain} 
\newtheorem{theorem}{\indent\sc Theorem}[section]

\newtheorem{proposition}[theorem]{\indent\sc Proposition}

\theoremstyle{definition} 

\newtheorem{remark}[theorem]{\indent\sc Remark}

\makeatletter
\def\address#1#2{\begingroup
\noindent\parbox[t]{7.8cm}{%
\small{\scshape\ignorespaces#1}\par\vskip1ex
\noindent\small{\itshape E-mail address}%
\/: #2\par\vskip4ex}\hfill%
\endgroup}%
\makeatother
\author{Pietro Corvaja,  Umberto Zannier}
\title{\uppercase{Examples of effectivity for integral points on certain curves of genus 2 }}
\date{\today}
\begin{document}

\maketitle



\section{Introduction}  This short article concerns a method to obtain effectivity for the search of integral points  on certain (sets of) curves of genus $2$.

As a preamble, it will be useful to mention that this draft is  intended merely to be a brief and very sketchy exposition of the said method, actually limiting the statements to a few examples. Not only we shall not give complete proofs, but often we shall omit even partial arguments. We plan to prepare in a near future a paper with more extensive statements and complete detail, also in collaboration with Davide Lombardo, who already performed several   verifications  used in this draft, and whom we thank.

However we hope that the present exposition at least will indicate the main lines of our method and the type of results that can be so obtained.

 \medskip

More precisely, we wish to illustrate just an example of application of a criterion of Bilu \cite{Bi}, to derive effectivity for integral points on certain sets of affine curves $X/\overline\Q$ of genus $2$, with smooth projective model $\tilde X$, so that $X=\tilde X-\{q_0\}$, $q_0$ being some (non-special)  point on $\tilde X$.  As said, future work will contain details and  more general applications (however of the same shape, in particular again with curves of genus $2$).

\medskip

The results will be far from giving effectivity for all such curves $X$. However they shall apply to a set of curves (and points $q_0$) which,   for instance, is  topologically complex-dense in the whole moduli space  for such data. (See Theorem \ref{T} below for a flavour of the statements that we can obtain. Also, see in particular the Notes at the end of \S\ref{SS.hints} and at the end of \S\ref{SS.concl}.)

 We shall construct morphisms   from these curves to $\G_{\rm m}^2$, with image  of increasing degrees.  We note that as the degree increases, we may say that the examples become `more interesting',  since they cannot be derived by substitution from a universal family. As a counterpart, there is the  negative feature in that  the relevant   curves will somewhat have increasing  fields of definition.

 On the other hand, to our knowledge this type of result did not previously appear in the literature. Actually,   the effective known cases of this type seem to be sporadic; for instance,   Serre's book \cite{Se}  explicitly  inquires about   possible effectiveness for curves  of this shape (see the remarks at the end of p. 116). In fact, this shape is somewhat significant, being  in a sense the simplest for which effectivity is not generally known.\footnote{Indeed these may be realized as  affine quartic curves of genus 2 with one point at infinity, whereas effectivity is known for all affine curves of genus $\le 1$, in particular for cubics.}

Further,  we do not know whether the mentioned   kind of  limitation is intrinsic, e.g. of the method. In this direction, we point out  that recent work of Landesman-Poonen \cite{LP} indeed suggests that Bilu's criterion  is subject to strong obstructions.\footnote{We note that the results in \cite{LP} refer to curves over $\C$. In principle by some kind of specialization pattern one is led to expect that the same should hold over $\overline\Q$; however some care is needed in taking for granted such a type of inference, as  is shown for instance  by remarkable examples like Belyi's Theorem.}

\bigskip

We recall that any curve  $\tilde X$ as above of genus $2$ may be obtained by adding two points at infinity to an  affine curve  defined by an equation $y^2=f(x)$, where $f$ is a polynomial of  degree $6$  without multiple roots.  This exhibits $\tilde X$ as a double cover of $\P_1$ ramified at $6$ points (called  {\it special points} or {\it Weierstrass points}). The datum of these points up to action of $PGL_2$ determines $\tilde X$ up to isomorphism.\footnote{Here and in the sequel we shall work with  $\overline\Q$ as ground field.}  (There are fairly complicated {\it Igusa invariants} for the resulting quotient space.) 

There is a natural (unique) involution on $\tilde X$ with quotient $\P_1$, defined by $(x,y)\mapsto (x,-y)$, and this fixes the six said  special points. 
 If $q_0$ is a special point then $X$ has a model in $\A^2$ given by an equation $y^2=f^*(u)$ where $f^*$ has degree $5$, and then  the integral points may be found effectively with well-known methods. \footnote{The criterion of Bilu's that we shall use may be read as including these cases, however with a construction of a cover which goes back to Siegel and differs from the  one that we shall introduce.}  The same effectivity is possible  for an affine  curve defined in $\A^2$ by $y^2=f(x)$, with $\deg f=6$ (actually any $f$), however in this last affine model  we are removing {\it two} points from  $\tilde X$, moreover related by the said involution.

Instead, our examples will involve a single non-special point $q_0$.

The (complex) algebraic curves of genus $2$ are parametrized by a quasi-projective variety denoted usually $M_2$. This is  a rational variety of dimension $3$. In the sequel the terminology `moduli space for curves of genus $2$' will refer to some quasi-projective variety in birational correspondence with $M_2$. So generically a point of such variety will represent a curve of genus $2$ and every curve of genus $2$, with the exception of a family at most of dimension $2$, will correspond to such a point. We shall also refer to $M_{2,1}$ as a quasi-projective variety  parametrizing  the pairs $(X,p)$ where $X$ is a curve of genus $2$ and $p$ is a point in $X$; it has dimension $4$.

\medskip

An instance of the results is represented by the following

\begin{theorem}\label{T}  Let $S$ be a finite cover of a dense   Z.-open subspace of a moduli space for curves of genus $2$. Let $\tilde X_s$ be the curve corresponding to $s\in S$. Suppose we are given  a prescribed (non-special) point $q_s\in \tilde X_s$  varying rationally as $s\in S$. Then there is a complex-dense subset $\Sigma\subset S$ of algebraic points on $S$ such that  for each $s\in\Sigma$ the integral points on $\tilde X_s-q_s$ are  effectively computable, over EVERY  number field.
\end{theorem}

NOTE: (i) The set $\Sigma$ itself  is `effective', so to say, for instance  in the sense that we may compute some element  of $\Sigma$ in any prescribed  open disk in $S$.  And we may establish whether a given algebraic point of $S$ belongs to $\Sigma$. The set $\Sigma$  is described by denumerably many explicit algebraic equations. The nature of these equations will be clearer from the arguments. We may say that the points of $\Sigma$ correspond to points of $S$ where a certain section of an abelian scheme  on  (a finite cover $S'$ of)  $S$, assumes values which are torsion in the appropriate fiber.

(ii) We can see $S$ as a finite-degree cover of an open affine  subvariety of $M_{2,1}$ whose projection to $M_2$ is of finite degree.

\medskip

We again stress that the present draft will not contain a fully detailed proof of this theorem; however we shall give many hints on it. A more complete treatment will follow in \S \ref{S.pencil}, where we consider only a pencil of curves inside the moduli space $S$. There we shall also indicate other types of results which are in the range of this method.

Before explaining the principles of our method, we pause to give some   general  considerations, to which we shall refer  in the sequel.

\subsection{Some explicit examples}

 Each of the mentioned  curves may be also seen as a quartic in $\A^2$, with a single (smooth) point at infinity and precisely one singular point (see below for a very simple argument). This shape phrases  the problem in classical diophantine terms. Let us see some examples of such quartics.

Consider the family, of genus $2$,  given by equations
\begin{equation*}
y^4+ay^2-xy-x^3+bx^2=0,
\end{equation*}
with parameters $a,b$.  This is found to be (generically) singular (only) at the origin, and to have a single (smooth) point at infinity, in a closure in $\P_2$. The genus is $2$ by the formula $g=(4-1)(4-2)/2$ for a smooth quartic, but the singular point subtracts $1$. 

We will treat this family in the mentioned forthcoming paper, performing the necessary (computational) verifications so as  to check that it satisfies the same conclusion as in Theorem \ref{T}.  But here we limit ourselves to point out a subfamily capable of effective treatment (and another problem in diophantine effectivity).

\medskip

AN EFFECTIVE CASE: Once in a certain expository paper (by Masser-Zannier), seeking  maximum simplicity,  the example $y^4-axy-x^3=0$ was proposed. But the authors missed to note that the existence of an automorphism of order $5$ ($x\to\theta x, y\to \theta^2y$) leads to  the affine model $y^5=w(w+a)^2$ where $w=x^2/y$. This last of course falls into the usual treatments. 
(Bogdan Grechuk  pointed out this case to Masser and Zannier, with a different and much longer  series of substitutions.) 

This very peculiar feature will not happen for the examples below (and it does not happen for the general equation displayed above), as is not difficult to check.

\medskip

ABOUT TRINOMIALS. The case of trinomial equations appearing in the last paragraphs, actually is always effectively solvable. Over the usual ring of  integers $\Z$, this was observed and proved by B. Grechuk, T. Grechuk, A. Wilcox \cite{GGW}, relying on Runge's theorem, on Baker's results on superelliptic equations, after a rather long argument with comparison of factorisations, and distinction into cases. In particular, due to the use of Runge, their proof, as it stands,  does not work over general number fields but only over $\Z$ (and even the rest of the argument does not take nontrivial units into account).  

Such proof in practice exploits the Newton polygon, so as to express $x,y$ multiplicatively  in terms of two other functions. A  shorter, different,  and completely general, proof goes as follows. Let $x^n+ax^ry^s+by^m=0$ be the equation, defining a curve $Z$, where the exponents are $\ge 0$ and $a,b\in k^*$. Put $\Delta=|nm-rm-sn|$. If $\Delta=0$ then  $u:=x^n/y^m$ and $v:=x^r/y^{m-s}$ are multiplicatively dependent and the equation $u+av+b=0$ shows that the curve $Z$   is the union of translates of tori in $\G_m^2$. If $\Delta\neq 0$,  the curve is irreducible, as can be seen by the substitution $x\mapsto t^m x, y\mapsto t^ny$ which transforms the equation into a binomial in $t$.  \footnote{ 
See Schinzel's book \cite{S} for a complete and subtle  theory of reducibility of trinomials over function fields, not needed in the simple case considered here.}  There is an abelian-group action on  the function field $\bar k(Z)$  of $Z$, defined as follows. Let $\theta,\eta$ be roots of unity of order (dividing) $\Delta$  such that $\theta^n=\theta^r\eta^s=\eta^m$; then we have an action of the group of such pairs $(\theta,\eta)$, defined by $(x,y)\mapsto (\theta x,\eta y)$, well-defined on $\bar k(Z)$.  
If $\theta=\zeta^p,\eta=\zeta^q$ for a fixed primitive $\Delta$-th root of unity $\zeta$, then this corresponds to $pn-qm\equiv pr-q(m-s)\equiv 0\pmod\Delta$.  Hence the group corresponds to integer vectors $(p,q)$ satisfying these congruences, which is a subgroup of $\Z^2$ of index $\Delta$, modulo $\Delta \Z^2$. The group leaves invariant 
  the functions $u,v$. Actually, the field of invariants is generated by $u$ and  $v$ (because the lattice generated by $(n,m)$, $(r,m-s)$ has index $\Delta$ in $\Z^2$).  Since $u,v$ are linearly related,  the field of invariants corresponds to a curve of genus $0$. Also, the group action induces an action on the points of a smooth complete model $Z'$ of $Z$, which leaves invariant the set of poles of $x$. Hence the affine curve for which we seek the integral points may be seen as  a fiber product of two cyclic covers of the affine line (possibly deprived of some points), and thus the known methods allow to determine the integral points, if their set is finite, or to parametrize them. It turns out that there may be infinitely many ones only if the group above is cyclic (which happens if an only if $\gcd(m,n,r,s)=1$). \footnote{The case of general Galois covers of the affine line is less direct than the special case appearing here. It still may be treated effectively, as in  independent results by Bilu \cite{Bi}  and Dvornicich-Zannier \cite{DZ}. Both of these results in fact appeared in previous publications, very difficult to trace back nowadays. Another independent argument had been given by Poulakis. Serre interpreted all of these arguments within a short proof, in a letter to  Bertrand of 4 June 1993. 
  See also \cite{Z2} for full and self-contained detail  for   the general case of genus $0$.}

\medskip

HYPERELLIPTIC MODEL: Intersecting the quartic with lines through the singular point one finds the hyperelliptic genus $2$ involution:
\begin{equation*}
(x,y)\mapsto (2{x^4\over y^4}-x, 2{x^3\over y^3}-y),
\end{equation*}
whence we have the usual sextic model
\begin{equation*}
\mu^2=a\lambda^6+\lambda^5+b\lambda^4-1,\qquad \lambda={y\over x},\quad \mu=\lambda^{-4}-x.
\end{equation*}

Of course this is to be interpreted as defining the function field of the relevant curves; indeed, the affine model with such equation has two points at infinity.

\subsection{\bf Perfect-powers values of algebraic functions}

We pause to point out a connection of the problem of effective integral points on such curves of genus $2$, with another natural problem in diophantine theory.

\medskip

PROBLEM.  {\it  Let $\pi: Y\to \A^1$ be  a (possibly ramified) cover of (affine) curves over a number field $k$, and suppose we are given  a rational function $f\in k(Y)$. The problem is to find  effectively  the integral (or rational) points $m$ on $\A^1$  such that they lift to an algebraic  point $\xi=\xi_m\in Y(\bar k)$  with the property that $f(\xi)$ is a perfect $p$-th power in $k(\xi)$. }

\medskip

This integer $p\ge 2$ may be fixed or variable, and the function $f$ is supposed not to be a perfect power in $k(Y)$.

\medskip

NOTE: Let $Y=\A^1$ and $\pi$ be the identity. Now  $f\in k(x)$,  and we have $\xi=m\in  k$,  so we are seeking the $p$-th power values (over $k$) of a usual rational function. These integral points (i.e. such values at integers of $k$) may be found  effectively (by Baker, as applied by Schinzel-Tijdeman), even letting $p$ be variable.

\medskip

In general, we can take the norm of $f$ down to $k(x)$ and obtain a rational function $F=F(x)$ on the affine line. Then we can often  apply  Schinzel-Tijdeman (see for instance the paper by Berczes, Evertse and Gyory \cite{BerGyo}). But a problem arises if $F$ is itself identically a perfect power (for instance a constant)  even if $f$ is not.

For any  fixed $p$ the issue may be dealt with through Siegel's theorem. In fact, if $f$ is not (identically) a $p$-th power  of a rational function on $Y$, the whole  graph of fields (going to a normal closure of $k(Y)(f^{1/p})/k(x)$) degenerates under specialisation at $\xi_m$, and we have an integral point on a certain curve (a quotient of the normal closure  by a suitable decomposition group), to which Siegel's theorem may be applied, and we can decide effectively  whether there are infinitely many integral points.\footnote{It seems that for $p>2$ we always have finiteness. Certainly for $p=2$ there are cases with infinitely many integral points, which can be found through equations of Pell type.} 
But, motivated by the positive results for powers, especially by the one  of Schinzel-Tijdeman,  we still may ask:  

QUESTION: is this problem completely effective ? Can we find these integral points when their set is finite ?

\medskip

CUBIC ROOTS. There is a link of (a special case of) this problem (case $p=3$) with the above problem of curves of genus $2$.  Indeed, we may see the above equation as a cubic in $x$ and solve with Cardan's formulae. For instance take $b=0$ in the above family so we have an equivalent  equation 
\begin{equation*}
x^3+xy+y^4+ay^2=0.
\end{equation*}
The formulae give, up to a constant,  cubic roots of $-y^4-ay^2\pm\sqrt{(y^4+ay^2)^2+(4/27)y^3}$ so of shape 
\begin{equation*}
\root 3\of {-y^4-ay^2\pm y\sqrt{(y^3+ay)^2+(4/27)y}}.
\end{equation*}

Corresponding to an integral point $x=r,y=s$ on the original curve, the  expression for which we have to take a cubic root  would produce a value in $\Q(\sqrt{(s^3+as)^2+(4/27)s})$.

Then, the issue is that we do not know, and is not necessarily true,  that the square root $\sqrt{(s^3+as)^2+(4/27)s}$ produces  a rational value, i.e. that  $(s^3+as)^2+(4/27)s$ is a square in $\Q$. Indeed, this is a stronger constraint and we could effectively find the integers $s$ for which this is so. 

Note that the norm of the relevant  function is a perfect cube, so we fall in the exceptional cases said above, where the known theory does not provide an effective answer.  These formulae provide an explicit link of this problem with effectivity for integral points on a curve of genus $2$.

\section{Our method - Curves of genus $2$}

We shall use a criterion of Bilu \cite{Bi}: \footnote{We apologize for the clash of the notation  $S$ both for our base space and a finite set of places; however the latter notation is so common that we prefer to forget about this, which is not likely to cause any confusion.} 

\medskip

{\bf Bilu's criterion}:  {\it If the group of regular morphisms $X\to\G_m$ has rank $\ge 2$ then finding the $S$-integral points on $X$ is effective   over any number field and any finite set $S$ of places}.

\medskip

Further phrasing: 

{\it If the group generated in $J_X$ by the differences of points at infinity (i.e. points in  $\tilde X-X$) has  co-rank $\ge 2$,   then there is effectivity for integral points}.

\medskip

There are still other equivalent versions.

\medskip

This is based on the following result:

{\it The equation $f(x,y)=0$ in $S$-units $x,y$ is effective}.

Proof of this: linear forms in logs + Puiseux series (see Bilu's papers or  \cite{BG}, Ch. 5,  for such an argument avoiding even an appeal to Puiseux theorem).

\medskip

As Bilu did in some cases (e.g. for modular curves, \cite{Bi}, \cite{BiI}), one may apply the criterion after taking a cover unramified outside the points at infinity (using the Chevalley-Weil Theorem). Then, as formulated explicitly by Bilu,  a  sufficient condition for effectivity becomes:

\bigskip

{\bf Bilu's improved statement}: {\it If $X$ has an unramified  cover $Y$  such that the differences of the lifted points at infinity  in $\tilde Y$ have co-rank $\ge 2$ in $Jac(Y)$, then $X(\O_S)$ is effectively computable (and actually so is $Y(\O_S)$).}

\bigskip

SIMPLEST CASE OF OUR METHOD

As above, the curve $\tilde X$ may be realized as a smooth complete model with function field given by $v^2=f(u)$ where $f$ has degree $6$. 

We define $X=\tilde X-\{q_0\}$ for a (non-special) point $q_0\in\tilde X$, to be chosen.  

We let $J_X$ be the Jacobian of $\tilde X$.

\medskip

We then construct an unramified cover $\pi:\tilde Y\to \tilde X$ cyclic of degree $3$, with points $p_1,p_2,p_3$ above $\infty=q_0\in X$.\footnote{For Bilu's criterion it would suffice that $\pi$ is unramified except above $q_0$. However we forget about this freedom in the present approach. The paper \cite{LP} shows that this freedom `generically' is not a significant improvement.}  

If the Galois group of $\pi$ is generated by $g$, of order $3$, we may label the points modulo $3$ so that $p_{i+1}=g(p_i)=:p_i^g$. 

 This cover $\tilde Y$ corresponds to taking a cube root of a rational function on $\tilde X$ \footnote{We suppose that the ground field contains a cubic root of $1$, actually let us take freely a finite extension of the ground field in the whole approach.}; considering its divisor and the fact that the cover is unramified, a simple analysis shows that this function corresponds to a point of order $3$ in $J_X$, and there are $80$ (though half of them yield the same covers). We omit here the relevant verifications, which will appear in a joint paper also with Lombardo.

\medskip

{\bf Applying Bilu's criterion}: In order to apply Bilu's criterion, we want  that  the $p_i-p_j$ are torsion in $J_Y$.  

\medskip

NOTE: A priori this is too difficult to achieve; in fact, $J_Y$ has dimension $4$, whereas the whole space of curves of genus $2$ has dimension $3$. The point at infinity yields one more freedom, but still a priori not sufficient to achieve (many)  examples. \footnote{In fact, see the paper \cite{LP} for grasping the limited hopes of such search.}

However there are three peculiar features in our  situation, which give more freedom, which we now  illustrate.

\medskip

{\bf A}.  First of all,  $J_Y\approx J_X\times J'$ (i.e. up to isogeny), since there is a surjective homomorphsm $J_Y\to J_X$, so $J'$ may be taken as the kernel.  

Note that $\pi\circ g=\pi$, so  the differences $p_i-p_j$ map to $0$ on $J_X$, so their classes in $J_Y$ (indicated with the same letter here) belong to  $J'$, which has dimension $2$. This is already a big gain. 

\medskip

{\bf B}. These differences  (as classes in $J_Y$)   are  found to be not generically torsion; note  however, that if they were, we could apply Bilu's criterion at once to the whole set of curves.  

Using the cyclic structure  it follows that  $p_3-p_1$ and $p_2-p_1$  are generically independent  over $\Z$ (i.e. as  $s$ and the curve $\tilde X_s$ vary).   Indeed any linear equivalence relation $ap_1+bp_2+cp_3\sim 0$, $a,b,c\in\Z$ with $a+b+c=0$,  implies $ap_2+bp_3+cp_1\sim ap_3+bp_1+cp_2\sim 0$, and it follows that $ac(p_2-p_1)\sim ab(p_1-p_2)\sim 0$, so the differences would be torsion.

However they are linearly dependent over $End(J_Y)$; indeed, $J_Y$ too has  an automorphism $g$ of order $3$ coming from the one of $\tilde Y$, and $\{p_1-p_2,p_2-p_3,p_3-p_1\}$ is  an orbit, so $p_3-p_2$ is the image $g(p_2-p_1)$ of $p_2-p_1$ (for  the previous labelling of the $p_i$). 

Hence if one of the differences  is torsion, any other one is also torsion. Therefore it suffices to make one of them torsion to achieve Bilu's conditions (in fact  $p_2-p_1, p_3-p_1$ are plainly independent over $\Z$ as divisors). This is a further gain and  (a priori)   now we have 4  or 3  freedoms (that is, depending on whether we vary $q_0$ or not) against only 2  constraints (i.e., the dimension of $J'$).

\medskip

We go ahead by describing a further pleasant  feature, which will be exploited in \S \ref{S.pencil} below. 

\medskip

{\bf C}. We gain even more simplicity and  freedom by the following reason.

It is not difficult to see that the hyperelliptic involution on $\tilde X$ lifts to $\tilde Y$.  
Actually it lifts in three ways, conjugate under the Galois group of the cover.  Then $\tilde Y$ becomes a Galois (ramified) cover of $\P_1$ with Galois group $S_3$. The curve $\tilde X$ is the  quotient of $\tilde Y$ by the subgroup of order $3$ in $S_3$.

Any of the three quotients of $\tilde Y$ by a transposition is found to be a curve of genus $1$, and these curves are isomorphic because the three transpositions are conjugate in $S_3$. This yields three morphisms from $J'$ to the cube of an elliptic curve $E$; if $\phi$ is one of them, the other ones are $\phi\circ g, \phi\circ g^2$ for a $g$ of order $3$. The sum of the three maps  from $J'$ is  easily seen to be zero, so we  remain with a morphism of $J'$ to a square $E^2$. A little elementary Galois cohomology (by the way for $S_3$) yields that $J'\approx E^2$ (i.e. up to isogeny). This can be seen also by checking that two of the maps are independent; indeed if they were not, then acting with $S_3$ we would find that each map is constant, which is absurd.  


--------------

NOTE: Not all these curves (that is, for other choices of the point of order $3$)  are isomorphic, in fact their $j$-invariants  are found to assume   $40$ values.  We owe these calculations to Lombardo.

QUESTION: are they all isogenous ?  

We haven't checked this. Maybe it is easy to see an answer.  In any case, this does not matter for our application.

NOTE: After this note was written Lombardo has found a negative answer (details to appear in the said forthcoming paper). He also showed us that in his joint paper with Garcia, Rizenthaler and Sijsling \cite{Lomb} the  decomposition, up to isogeny, of the jacobian of $\tilde{Y}$ as $J_X\times E^2$ appears as a special case of their work on decompositions of jacobians via Galois covers (see in particular Table 7 therein).
--------------

\medskip

This isogeny decomposition $J_Y\approx J_X\times E^2$, hence $J'\approx E^2$,  yields further freedom, since we may seek for instance a section which is torsion identically on one elliptic factor, or such that its two elliptic components are identically linearly dependent, and then we have only to prescribe torsion on one of the two   factors, which is a weaker constraint, allowing to lower the dimension of the base.
Indeed we shall later  give an example of a 1-dimensional family, based on exactly this type of  construction. 

\medskip


\subsection{Hints for a  proof of Theorem \ref{T}} \label{SS.hints} 

Let $S$ be as in the statement, and let us be given a non-special algebraic section $s\mapsto q_s$. 
So for each smooth complete curve $\tilde X$ of genus $2$ (except possibly for the curves corresponding to finitely many surfaces  in $M_2$), we have a  finite non-empty set (with uniformly bounded cardinality) of  $s\in S$ such that $\tilde X_s$ is isomorphic to $\tilde X$, and we have a    (generically non-special) point $q_s\in \tilde X_s$, where $q_s$ is given with coordinates which are rational functions on $S$.   (For our purposes here, it suffices   to use such naive notions of moduli space.)

\medskip

For $s\in S$ we may construct as above a cover $\tilde Y_s$ of $\tilde X_s$, after a choice of a point of order $3$ on $J_{X_s}$. This $\tilde Y_s$ in fact which will be defined over a  suitable finite cover of $S$ (depending for instance on the choice of  such point of order $3$). \footnote{Explicit formulae can be computed, and part of them will appear below. Lombardo has performed complete computations, which will appear in a joint paper with the present authors.} The Jacobians of the curves $\tilde Y_s$ determine an abelian scheme over this cover of $S$. 

Replacing $S$ by a  further finite cover (and denoting it again by $S$ for convenience), we obtain two sections from $S$  to this abelian scheme,  i.e. the classes of the differences $p_2-p_1$, $p_3-p_1$ (where $p_i=p_{is}$), viewed as points of $J_{Y_s}$.  As noted above, these sections are dependent over $End(J_{Y_s})$. If they would be dependent over $\Z$ we would have a nontrivial  linear equivalence relation $a_1p_1+a_2p_2+a_3p_3\sim0$ with $a_1+a_2+a_3=0$. Acting with the group $<g>$ the $a_i$ get cyclically permuted and we eventually find that  the class $p_2-p_1$ is identically torsion on $S$.  
If this were  the case then we could apply Bilu criterion to all members of the family, leading to a stronger theorem. So let us suppose this is not the case.\footnote{It is however possible to exclude that  this  happens on the whole $S$, no matter the choice of the section. It follows from a general theorem in \cite{CMZ} that this is true  except for a finite number of sections $q_s$. 
However below in \S \ref{SS.nt} we shall sketch a  proof of the stronger  assertion.} 

\medskip

Now, even after we have prescribed the  family of points at infinity $q_s$ (depending on $\tilde X$),   losing in this way one freedom, we are left with  a  base  space $S$ of dimension $3$.
Then by the theory of the Betti map  we can make one of the sections torsion in a dense set $\Sigma$, so the other section will be torsion as well by remark (B) above.  The former kind of verification may be generally not easy to achieve (see e.g. the paper \cite{ACZ}), but here the situation is simpler, since the relevant abelian varieties are (isogenous to) squares of elliptic curves, by comment (C) just above. This case is easier and treated in particular in \cite{CMZ}. 
But we need not even appeal to that paper, the main points of the argument being as follows:

\medskip

\noindent{\bf Sketch of the argument with the Betti map}:  Suppose we have a doubly elliptic scheme $E^2/Z$ without constant part (i.e., the curve $E$ has not constant $j$-invariant as a function on $Z$), over a(n affine)  surface  $Z$, with a section $\sigma$. Then we have to prove that  the points $z\in Z$ when $\sigma(z)$ is torsion are complex-dense in $S$. Now, writing $\sigma=(\sigma_1,\sigma_2)$, by the elliptic case (which is not difficult  and treated e.g. in \cite{Z}, see Notes to Ch. 3), there are a complex-dense set of algebraic curves in $Z$ such that $\sigma_1$ is identically torsion on each of these curves.  (In fact, it suffices first to prove this on a `general' curve $C$ in $Z$, which is a known case. Then each point of $C$ where $\sigma_1$ is torsion will yield a full algebraic curve in $Z$ on which $\sigma_1$ will be torsion of the same order.)  In turn, for each of these curves, by the same principle, there is a complex-dense set of points where $\sigma_2$ is torsion, concluding the argument.

\medskip

NOTES. (i) The torsion orders so obtained will grow to infinity, unless the section  is  identically torsion (on the whole base $S$); this last possibility 
shall be excluded  in \S\ref{SS.nt} below (see  Proposition \ref{P.nt}). 

(ii) If the curve $E$ is isotrivial, so may be assumed constant,  the assertion is not generally true but remains true (with a very easy proof) provided only that the section is not (a non-torsion) constant. That this case cannot happen in our case follows from the fact that the elliptic family coming from the construction has not constant $j$-invariant  on the whole base  $M_2$ (as shown by an explicit computation of Lombardo, appearing in a future joint paper and not inserted here). See also  the final considerations at the end of this article.   For the   purposes of the construction in \S\ref{S.pencil} below, involving a pencil,  this is immaterial since the section   will shown explicitly not to be constant. (See the step S8 in the construction below, for this specific issue.) 

\medskip


So the set $\Sigma$  mentioned in the theorem in practice  is the set of points of  $S$ where the differences $p_i-p_j$ become torsion in $Jac(Y_s)$.

\section{A pencil with infinitely many effective cases}\label{S.pencil}

We give an explicit construction by a direct geometrical process. We start, for instance,  with the Fermat  curve of genus $1$ given by ${\mathcal F}:u^3+v^3=1$. In this pencil the relevant elliptic curves will be all isomorphic to ${\mathcal F}$.

We  proceed according to  the following steps.

\medskip

S1. We project $E$ from a variable point  $(0:\alpha:1)$
 (lying on the line $T=0$ at infinity in homogeneous coordinates $(T:U:V)$)  to the line $v=0$, i.e. the $u$-axis.
 
\bigskip
 
 S2. There are $6$ branch points for this projection (by Hurwitz), and a calculation shows that  they are given by $v=0$ where 
 \begin{equation*}
   v^2=(\alpha^3+1)u^6-2u^3+1, 
 \end{equation*}
is the equation of a curve of genus $2$, which will define our $\tilde X_\alpha$ (where we disregard $\alpha^3=-1$). At the moment  we have a scheme over the affine $\alpha$-line.

 It may be shown that there are no other automorphisms  of $\tilde X$ than the obvious ones (we mean generically in $\alpha$, for instance there is a nontrivial one for $\alpha=0$),  so that we do not fall neither  into a Thue nor in a  hyperelliptic case when we take $X$ as $\tilde X$  minus a single non-special  point.  Indeed, the issue of automorphism is a relevant one: if for instance there is an automorphism fixing the point at infinity, then the quotient has genus $0$ or $1$, and our curve is a cyclic cover of it, and hence there is an effective algorithm by known facts.

 It is not too difficult to see as well that the curve-family  is non-constant, that is,  the curves in the family are generically non-isomorphic (not even after extension of the base).   We omit these verifications in this preliminary draft,  which seem however  not entirely evident and requiring at any rate some checking. 

\bigskip

S3. We have $(\alpha^3+1)u^6-2u^3+1=P^2-Q^3$ with $P=u^3-1$, $Q=-\alpha u^2$.  

We use this representation to define a cubic cyclic unramified cover of $\tilde X_\alpha$. In general, such a cover  can be obtained  by choosing a point of order $3$ on $J_{X_\alpha}$; we omit such construction here and only give the results of the computation (which, as said,  will appear elsewhere). In the present case this cover  will  continue to be defined on  the $\alpha$-line. Indeed, it may be shown that the said choice of the point of order $3$ corresponds to a choice of the polynomials $P,Q$ such that $(\alpha^3+1)u^6-2u^3+1=P^2(u)-Q^3(u)$.  (To be precise, in fact $P,Q$ yields a bit of more information, immaterial for the present purposes, we shall make this precise in the mentioned future paper.)  In our case this choice is already implicit in the opening construction; this explains that the field of definition does not presently increase. 

More precisely, omitting the subscript $\alpha$ for notational convenience, the function field of the cover $\tilde Y$  may be  obtained by adding to the function field $k(u,v)$ of $\tilde X$ a cube root of the  function $v+P(u)=v+u^3-1$, which we denote now $w$ so $w^3=v+u^3-1$, and the function field of $\tilde Y$ is generated by $u,w$ over $\overline\Q(\alpha)$. 

\medskip

Let $g$ generate the group of automorphisms of $\tilde Y$ as a cubic cyclic cover of $\tilde X$. Then $g$ may be chosen so to act as $(u,w)^g:=g(u,w)=(u,\theta w)$ where $\theta$ is a given primitive cubic root of $1$.

\bigskip

S4.   We can then apply suitable calculations (again, omitted here) for obtaining the   curves of genus $1$ mentioned above,  as quotients of $\tilde Y_\alpha$ by an involution lifted from the canonical one on $\tilde X_\alpha$.
 In this case we find the cubic model of one such   curve  given by
 \begin{equation*}
E_\alpha:\qquad  z^3-3Q(u)z-2P(u)=z^3+3\alpha zu^2-2(u^3-1)=0.
 \end{equation*}
We have used affine coordinates but for $E_\alpha$ we implicitly mean a projective closure of this cubic in $\P_2$, whose equation  reads $Z^3+3\alpha ZU^2-2U^3+2T^3=0$,   in projective coordinates $(Z:U:T)$. The discriminant is found to be $\alpha^3+1$ up to a constant.  This $E_\alpha$ is isomorphic (over some extension of $\Q(\alpha)$)  to the Fermat cubic chosen at the beginning.\footnote{The shape  of these equations is cubic form$=$ constant. By sending through $PGL_2$  the three zeros of the cubic  form to the zeros of a fixed one of them, one obtains explicitly a projective isomorphism with a constant curve.} 

These $E_\alpha$ are only curves of genus $1$ for the moment, since we have not yet chosen an origin of $E_\alpha$, but we shall do that very soon.

\medskip

Putting $z:=w+Q(u)/w=w-\alpha u^2/w$, we have a map $\phi=\phi_\alpha:\tilde Y_\alpha\to E_\alpha$  given by $\phi(u,w)= (u,z)$, as is easy to check. 

\medskip

Such  construction gives a pencil $\tilde Y_\alpha$ of  curves of genus $4$, with $\tilde Y_\alpha$ a(n unramified) cubic  cyclic cover of $\tilde X_\alpha$, and a morphism  $\tilde Y_\alpha\to E_\alpha$ (the curves $E_\alpha$ being all isomorphic to $\mathcal F$). 

In fact, another morphism to $E_\alpha$ is simply obtained by composing with $g$, so we have a morphism $\sigma: \tilde Y_\alpha\to E_\alpha^2$ given by 
\begin{equation*}
\sigma(p)=(\phi(p),\phi(p^g)),
\end{equation*}
where $p\in \tilde Y_\alpha$, denoting $p^g:=g(p)$, as usual.

Note that we can also use $g^2$ in place of $g$, but note that the divisors on $E_\alpha$ given by $\phi(p)+\phi^(p^g)+\phi(p^{g^2})$ are all linearly equivalent, since it is easy to check that they come from $\P_1$. Hence this sum-map is constant for every choice of an origin and acting with $g^2$ does not give anything new.

The map $\sigma$ induces a map from $Jac(\tilde Y_{\alpha})$ to $Jac(E_\alpha)$, which is surjective; in fact the involution on $\tilde X_\alpha$ may be lifted in three ways to $\tilde Y_\alpha$, which are conjugate under $g$, and this yields the maps  from the Jacobian (as has been said in the discussion (C) above).

\bigskip

S5. On this  family $\tilde X_\alpha$  we now choose our section $q_\alpha\in \tilde X_\alpha$ giving the point at infinity (to be removed so as to obtain the affine curve  $X_\alpha$).  We do this as follows.

Let $p_1,p_2,p_3$ be the fiber above $q_\alpha$ in $\tilde Y_\alpha$, where  $p_i=p_{i\alpha}$,  the labelling being  modulo $3$ being  chosen so that  $p_{i+1}=p_i^g$. 

We prescribe that $p_2=p_{2,\alpha}$ is a flex on $E_\alpha$, which will be chosen as origin. For a generic $\alpha$, the $u$-coordinates of the flexes are found to be defined by the equation
 \begin{equation*}
 (\lambda^3+3\alpha\lambda-2)u^3=-2\quad \hbox{where \quad $\alpha\lambda^2-2\lambda-\alpha^2=0$.}\footnote{The Hessian matrix has rows $(6z,6\alpha u,0)$, $(6\alpha u, 6\alpha z-12u,0)$, $(0,0,*)$.}
  \end{equation*}  There are also three flexes at infinity (i.e. with $T=0$). Note that for $\alpha=0$ the flexes are instead the points with $zu=0$ plus the three points at infinity.

 A flex at a finite point is defined either over $\Q(2^{1/3})$,  if $\alpha=0$, or over a cubic extension of $\Q(\sqrt{\alpha^3+1})$, that may be made explicit by the equations just given. This yields a cover $B$ of the $\alpha$-line over which $q_\alpha$ and the $p_i$ are   defined. 
 
 {\bf SMALL WARNING}: So in fact we should not use $\alpha$ as a subscript, but, since the beginning,  a point $\beta$ of $B$; however we omit this precision here.
 
 \medskip

 Note that by the above formulae, if $u_\alpha$ is the $u$-coordinate of the point $q_\alpha$, then $p_1,p_2,p_3$ will be sent to  the points on $E_\alpha$ with $u$-coordinate $u_\alpha$, and they will be distinct as soon as $\alpha$ is not a zero of the discriminant, i.e. $\alpha^3\neq -1$. 
 
 So, generically,  with the actual choice the point $q_\alpha$ will not be a special point on $\tilde X_\alpha$.

\bigskip

S6. According to the above settings and equations, the section $p_2-p_1$ on $Jac(\tilde Y_\alpha)$ is sent by $\sigma$ to the point of $E_\alpha^2$ given by 
\begin{equation*}
\sigma(p_2)-\sigma(p_1)=(\phi(p_2)-\phi(p_1), \phi(p_3)-\phi(p_2)).
\end{equation*}

Now,  let us note that the two components are equal: in fact, this amounts to $\phi(p_1)+\phi(p_2)+\phi(p_3)\sim 3\phi(p_2)$ on $E_\alpha$ (or the corresponding equality with $\phi(p_2)$ as origin, but that is immaterial for this).  In turn, this linear equivalence (or equation) holds because $\phi(p_2)$ is a flex, and the $\phi(p_i)$ are the intersections of $E_\alpha$ with  the line $u=u_\alpha$.

\bigskip

S7. We now want to check that the section is not identically torsion. If this were the case, then (acting with $g$) all the differences $\phi(p_i)-\phi(p_j)$ would be torsion.  

Now let us specialize to $\alpha=0$, which is a point of good reduction for the elliptic family, indeed $E_\alpha$ reduces to $Z^3-2U^3+2T^3=0$, which is a twist of the Fermat curve.\footnote{Note that $\alpha=0$ is not of good reduction for the family $\tilde X_\alpha$, but this is immaterial for the present purpose.}   For this curve the flexes are precisely the points with $ZUT=0$. The flex that we have chosen goes to a flex with $U=0$, and all $\phi(p_1),\phi(p_2),\phi(p_3)$ specialise to flexes. Then the differences have order $3$ in the reduced curve. By general reduction theory, if the section $\phi(p_2)-\phi(p_1)$ would be torsion, the order would divide $3$. But this is not the case, since generally $\phi(p_1)$ is not a flex, and is distinct from $\phi(p_2)$.

We conclude that the section is non-torsion.

\bigskip

S8. We now check that the section $\phi(p_2)-\phi(p_1)$ assumes torsion values for a complex-dense set of $\alpha$. If the elliptic family $E_\alpha$ would be non isotrivial, this would be a consequence of the theory of the Betti map, an  easy case presented e.g. in \cite{Z}, Notes to Ch. 3. But instead all curves $E_\alpha$ are isomorphic to the Fermat curve ${\mathcal F}:u^3+z^3=1$ and hence we have a(n iso-)constant family. The section then leads to a rational map from a dense Z.-open subset of a certain cover $\tilde B$ of $\P_1$ to $\mathcal F$; if this map is non-constant, then its image contains a cofinite set in $E$ and our assertion follows immediately. On the other hand if the section would be constant, it would be torsion because it is so for $\alpha=0$; but this contradicts what has been proved at the previous point, concluding the argument.\footnote{That the section to $\mathcal F$ cannot be identically constant follows also by direct calculation, on computing an explicit isomorphism of $E_\alpha$ with $\mathcal F$; since such an isomorphism may be found sending a  flex  to a flex, the argument at bottom is the same as specializing $\alpha=0$.} 

\bigskip

S9. Since by S6 above the two components of the section are equal, and since by S8 the first component assumes torsion values on a dense set, the section itself to $E_\alpha^2$ assumes torsion values on a complex-dense set in the relevant cover $B$ of the $\alpha$-line.

\subsection{Conclusion}\label{SS.concl}   Let us go ahead by using points $\beta\in B$ in place of $\alpha$, for more precision. We have a cover $B$ of the $\alpha$-line with a map $\delta:B\to\P_1$ where we denoted  $\delta(\beta)=\alpha$. 

Now we  have a family $X_\beta=\tilde X_\alpha - q_\beta$ as above. The construction gives a pencil $\tilde Y_\alpha$ of  curves of genus $4$, with $\tilde Y_\alpha$ a(n unramified)  cover of $\tilde X_\alpha$. Correspondingly we have an affine curve $Y_\beta:=\tilde Y_\alpha-\{p_1,p_2,p_3\}$  and a map $Y_\beta\to E_\alpha$ (the curves $E_\alpha$ being all isomorphic to $E$). 

The sections $p_2-p_1,p_3-p_1$ to $Jac(\tilde Y_\alpha)$ become both torsion in a complex-dense subset $\Sigma$ of (algebraic) points of $B$. For $\beta\in B$ Bilu's criterion applied to $Y_\beta$ yields effective finiteness for the $S$-integral points on $Y_\beta$ over any prescribed number field (and finite set $S$ of places). Since $Y_\beta$ is an unramified cover of $X_\beta$, the Chevalley-Weil theorem yields effective finiteness for the $S$-integral points on $X_\beta$ as well.

\medskip

NOTE: We may represent  $X_\beta$  as a pencil of quartics in $\A^2$, with a single singular point in $\A^2$ and a single point at infinity   in a smooth model $\tilde X_\beta$. We do not perform the computation here, but we reproduce an argument to show that this representation  is possible. (We do not call this an `embedding' since there is a singular point; however this is quite immaterial for our purposes since the singular point is not at infinity.) 

\medskip

======================

PAUSE: 

{\tt Representing an affine  curve of genus $2$ : $X=\tilde X-q_0$, $q_0$ non-special,   as a plane quartic with $q_0=\infty$}.

We simply consider the linear system $L(3q_0)$ (of rational functions on $\tilde X$ with at most a triple pole at $q_0$) which is generated by $1,z$ with $z$ a suitable rational function of degree $3$ (since $q_0$ is non-special). Then $L(4q_0)$ is generated by $1,z,w$ with a $w$ of degree $4$. Since $z,w$ have coprime degrees, they generate the function field of $X$. Since they have a unique pole they provide a map $X\to\A^2$ with a single point at infinity.  The extended map $(1:z:w)$ to $\P_2$ provides the projective plane model, with coordinate  functions   $1/w$ and $z/w$ at infinity, both vanishing therein; since the latter has a simple zero, and since these functions have no other common zero, the curve is smooth there. 

Let us show that the irreducible polynomial equation $f(z,w)=0$ relating $z,w$ is of total degree $4$. First of all, it has degree $4$ in $z$ and $3$ in $w$ (because $z,w$ have degrees resp. $3,4$ as functions on $\tilde X$).  Now, the space $L(15q_0)$ contains the $15$ terms $1,z,z^2,z^3,z^4,z^5, w, wz, wz^2, wz^3, w^2, w^2z, w^2z^2, w^3, w^3z$. By Riemann-Roch they must be linearly dependent, giving a nontrivial polynomial relation $g(z,w)=0$ where $g$ is a linear combination of the said terms. This $g$ must be divisible by $f$, proving that there is a relation without the term $z^5$, and of total degree therefore at most $4$ (and hence exactly $4$) as stated.\footnote{It seems that in  this case the usual  analysis using the spaces $L(mq_0)$, to embed the curve, needs the supplementary consideration involving $f,g$. Or else one can argue with $w+cz$, $c\in\C$,  in place of $w$ to deduce that $f$ has total degree $4$.}

This quartic model will  necessarily have  a single singular  (nodal) point in $\A^2$ by the genus formula. So in fact we have not truly an embedding of the smooth curve, but for the matter of integral points this is immaterial, since the singular point (which is always special)  is not at infinity.

======================

\medskip

{\bf Remark}. Note that the conclusions  are  achieved through a morphism from $Y_\beta$ to an irreducible  curve defined by an irreducible equation $F(u_\beta,v_\beta)=0$ in $\G_m^2$, not a translate of a torus, where $u_\beta,v_\beta$ are functions on $Y$ with divisors resp. $n(p_2-p_1)$, $n(p_3-p_1)$, where $p_i$ are the points of $Y$ above $q_0$.  Then these functions are multiplicatively independent modulo constants, and assume $S$-unit values at the integral points of $Y$ which are lifted by Chevalley-Weil from those on $X$. By Puiseux+linear forms this gives effective finiteness, as in the proof of the criterion of Bilu.

\medskip

 {\tt Degree of the image curve}. Let us check that for large torsion order of the relevant section, the degree of the image curve $Y^*:F=0$ also increases. Suppose then that $n$ is the exact (common) torsion order, so for every $m\ge 2$ the functions $u_\beta,v_\beta$ are not both perfect $m$-th powers in $k(Y)$, where in the following argument we omit the subscripts $\alpha,\beta$  for convenience. 
 
   We have an extension $k(Y)/k(u)$ of degree $n$, totally ramified below $p_1$ and $p_2$, and we have an intermediate extension $k(u,v)/k(u)$. The extension $k(Y)/k(u,v)$, say of degree a divisor $m$ of $n$,  is totally ramified below $p_1,p_2,p_3$ (at least).
   
   Suppose first that the genus of $Y^*$, i.e. of $k(u,v)$,  is $>1$. Then by the Hurwitz formula we would have $6=2g(Y)-2\ge 2m +3(m-1)=5m-3$, hence $m=1$, $k(u,v)=k(Y)$,  and $F$ would have degree at least $n$, hence exactly $n$.  
   
   Suppose now that $g(Y^*)=1$. Then again we obtain $6\ge 3m-3$, hence $m\le 3$ and we have a lower bound $n/3$ for the degree of $F$. \footnote{In  case of equality we would have a map of degree $3$ from $Y$ to an elliptic curve; then this curve would be isogenous to $E$, and several deductions would be possible. We do not go ahead, but maybe further  inspection could lead to the impossibility of this case.} 
   
   Suppose finally that $g(Y^*)=0$. Then $k(u,v)=k(z)$ for some function $z\in k(Y)$. After a homography we could assume that $u=z^q$ and $v=c_1(z-c_2)^q$ for  constants $c_1,c_2\neq 0$ and $qm=n$. The divisor of $z$ would satisfy $q\cdot \div(z)=\div(u)=n(p_2-p_1)$ so $\div(z)=(n/q)(p_2-p_1)$. Since we are assuming that  $n$ is the exact torsion order we conclude that $q=1$, and then $c_3u+c_4v=1$ for nonzero constants $c_3,c_4$. By the functional $abc$-theorem we deduce $n\le 2g(Y)+1=9$, so $n$ is bounded, concluding the argument.


\medskip

{\bf NOTE (a)}.  The integer $n$ may be taken as the common torsion orders of the two sections, which is also the torsion order of each of them, since they are dependent over $End(J_Y)$. 

It  grows to infinity, so it is not possible to  perform this procedure in a uniform (bounded-degree) way for all curves in the pencil. (In fact, if $n$ is the precise torsion order of  our  divisors in $J_Y$, the said morphism cannot be factored nontrivially  through a morphism of smaller degree to $\G_m^2$.) 

This is immaterial for the result itself, but it shows that one cannot  get this proof by universal equations through the family (as would be possible for instance for the affine curves $y^2=f(x)$ with two points at infinity). The situation somewhat reminds of Belyi theorem...

On the other hand, the same method is flexible in some aspects, and will give also examples of pencils where the relevant sections are identically torsion. For such  pencils we shall have effectivity for all algebraic points. However this effectivity will be in a sense less interesting, because it will correspond to a same type of diophantine equation (in the same way that all equations of type $y^2=f(x)$ may be effectively solved in integers).  These examples will appear in the mentioned paper in preparation with Lombardo.

\medskip

{\bf NOTE (b)}: The parameter will run on a curve of bounded degree over $\beta\in\P_1$. The relevant values of $\beta$ will be algebraic numbers with the following properties: 

 (i) of degree tending to $\infty$, 
 
 (ii) topologically complex-dense, and 
 
 (iii) of bounded height.

(iv) Moreover it will be effectively possible to check whether a given value of $\beta$ is or is not in our set.

\medskip


\subsection{Discussion about  the non-torsion of the section for the general family}\label{SS.nt}

We have seen in the hints of proof presented in \S \ref{SS.hints} that IF  a certain section of a doubly elliptic scheme, associated to the choice of the section $s\mapsto q_s\in \tilde X_s$ for points at infinity, is identically torsion, then we have  effectivity of integral points for the full set of affine smooth curves  $X_s$ of genus $2$ defined by $X_s:=\tilde X_s-\{q_s\}$. We want to prove now that this cannot happen, no matter the  choice of the section at infinity $q_s$.  This represents another limitation  of  the method. 


\begin{proposition}\label{P.nt} Notation being as in Theorem \ref{T}, there is no algebraic family $\{q_s\}_{s\in S}$, $q_s\in\tilde X_s$,  such that the doubly elliptic section  constructed as above is torsion on the whole $S$.
\end{proposition}

NOTE:  (i) As said, the paper \cite{LP} goes in a somewhat similar direction as this statement, namely providing limitations to Bilu's criterion. However there are substantial  differences with the present proposition. In a sense, the main result   of \cite{LP}  is much more general, as it concerns arbitrary covers of arbitrary curves of genus $>1$; and  it is different also because it starts with a {\it given} curve (rather than the whole moduli space) and proves that for a general choice of the point at infinity, Bilu's criterion will not work. Instead, here  we work only with cyclic cubic  covers, but we let the curves vary,  in genus $2$,  and we suppose that the choice of the point at infinity is given.

Also, the paper \cite{LP}  deals with objects defined over $\C$ (see a previous footnote for a comment on this).  
Here we may work over $\overline\Q$. Note that assuming that the value of the section at $s$ is torsion  for each $s\in S(\overline\Q)$, amounts to assume that it is torsion at a generic point: indeed, if the generic value is not torsion, known arguments of Silverman  (or see Demianenko-Manin in the isotrivial case) prove that the set of such $s$ has bounded height.

(ii) By appealing to a general result, we can easily prove that the Proposition holds up to   finitely many choices of $q_s$. For completeness (and possible use in other similar issues) we shall offer this short argument as well, before the proof of the stronger assertion.

\medskip

As a preamble,  suppose that a family as in the statement exists; this corresponds to having a  subvariety $S$ of $M_{2,1}$ projecting with a map of finite degree to $M_2$, so $S$ is a certain hypersurface in $M_{2,1}$. 

Also, the construction that we have described gives an abelian  scheme above a Zariski-open dense subset of a finite cover of $M_{2,1}$, with fibers which are squares of elliptic curves. Moreover, we have a section of this scheme, and its image is not contained in any group-subscheme (since the components of the section are  generically independent, as has been remarked several times). 

\medskip

====================

{\tt Pause: finiteness of possible exceptions}. 
An algebraic family as in the statement would correspond to a hypersurface where this section becomes torsion. However then we may apply Theorem 1.1 of \cite{CMZ}, which asserts the finiteness of such hypersurfaces, and this yields just the desired conclusion as in Note (ii) just above.

====================

\medskip

  \begin{proof} [Proof of Proposition] 
 Let us now prove the stronger assertion stated in the proposition, namely that such sections do not exist. In the notation above,   it suffices to prove that for every  choice of the  section $q_s$  at infinity  on $\tilde X_s$, the section $\phi(p_2)-\phi(p_1)$ to the relevant elliptic curve,  is not identically torsion; recall that  these  elliptic curves  $E_s$ are defined by an (affine) equation
\begin{equation}\label{E.cubic} 
z^3-3Q_s(u)z-2P_s(u)=0,
\end{equation}
where the curve $\tilde X_s$ of genus $2$ is defined (affinely) by $v^2=f_s(u)$, so that  $f_s(u)=P_s^2(u)-Q_s^3(u)$; here the polynomials $P_s,Q_s$ have  respective degrees $3,2$.

We remark that  any  polynomial $f$ of degree $6$ without multiple zeros is likewise representable, in finitely many ways, and the representations are in  a certain correspondence with the cyclic covers constructed above (again we omit such verifications here, but they will appear in the mentioned paper in progress). 

Note also that changing $f_s(u)$ with  
$(cu+d)^6f_s((au+b)/(cu+d))$ (where the coefficients depend on $s$, and $ad\neq bc$) yields an isomorphic curve. The same then happens changing $P_s,Q_s$ accordingly with the same transformation. 

By such a transformation (actually we need only a translation on $u$) we may  assume that the $u$-coordinate of $q_s$ is $0$. We then perform such a translation;  this is only to simplify the notation and to reduce the number of involved variables. (Of course assuming  that $q_s$ is non-special  amounts to $f_s(0)\neq 0$ (but this does not matter now.)

\medskip

It is relevant to recall that the points $\phi(p_i)$, $i=1,2,3$   (where we  again omit to indicate the dependence on $s$ for notational convenience), are the three points on $E_s$ with the same $u$-coordinate as the point at infinity $q_s$, i.e. the point  to be removed from $\tilde X_s$. Hence with the present normalization the points  $\phi(p_i)$ have $u=0$ and thus are defined by $z^3-3Q_s(0)z-2P_s(0)=0$, which has no repeated roots if $q_s$ is non-special.

\medskip

{\tt Going to a subspace with all isomorphic elliptic curves}. Let us now assume that the family is identically torsion on $S$; we have to prove that this is impossible.  

We can fiber the threefold $S$ with surfaces, denoted  $S_c$,  where the $j$-invariant of  $E_s$ is equal to $c$.  Let us restrict our attention to one (general enough) such $S_c$. Then for $s\in S_c$ the cubics defined by \eqref{E.cubic}, or rather their projective shape in homogeneous variables $(T:U:Z)$ (where $u=U/T, z=Z/T$),  are all projectively isomorphic, so any two of the resulting equations in $T,U,Z$ are related  
by a $3\times 3$ invertible linear substitution on the homogeneous variables. 

Let us now choose a flex on $E_s$ such that generically on $S_c$ its $U$-coordinate and $Z$-coordinates  are both  $\neq 0$. Note that such  flex exists (there are in fact $9$ flexes) and can be sent to $(0:1:1)$ by a projective transformation depending algebraically  
on $S_c$, represented by a matrix with rows of the shape $(*,*,0)$, $(0,*,0)$, $(0,0,*)$. Taking a finite cover of $S$ such that the flex becomes rational on the cover, the transformation too becomes rational; also,  composing with such a transformation does not destroy the  said property of $q_s$ (to have vanishing $U$-coordinate), and the polynomials $P_s,Q_s$ are changed according to an element of $PGL_2$, so the shape of the equation \eqref{E.cubic} is not changed. This further entails that the curve $\tilde X_s$ does not change its isomorphism class. Hence  we may assume that all curves $E_s$, $s\in S_c$, pass through $(0:1:1)$ having a flex there. We can also agree to choose that flex as an origin for all such curves, which then become elliptic.

\medskip

Now, let us fix (arbitrarily) one of the curves corresponding to a given  point $s_0\in S_c$, say $E^*:=E_{s_0}$. We have a projective transformation $L_s$ sending $E_s$ to $E^*$. The line $U=0$ is sent by $L_s$  to a line meeting $E^*$ in two  torsion points (and hence three). Therefore  by continuity this line must be constant as $s$ varies, so the $L_s$ all fix the line $U=0$. And the $L_s$ fix the flex $(0:1:1)$ as well, again by continuity (since the flex $(0:1:1)$ of $E_s$ is sent by $L_s$ to a flex of $E^*$, and there are $9$). 

Further, by composing (on the left) with a linear projective transformation with rows $(1,0,0), (0,\rho,0), (0,0,\rho)$, $\rho\neq 0$, does not destroy any of the shapes, fixes the flex, and fixes the isomorphism class of $\tilde X_s$ (it sends $P_s(u),Q_s(u)$ resp.  to $P_s(\rho u), Q_s(\rho u)$). We can now choose $\rho=\rho_s$ such that one intersection of $U=0$ with $E_s$ is constant in $s$, say equal to $(1:0:1)$.  (In fact, the roots of $z^3-3Q(0)z-2P(0)=0$ are finite and distinct.)

\medskip

Now, moving on a curve $S'_c$ inside $S_c$ we can achieve that also another intersection with the line $U=0$ is constant in $s\in S'_c$. But now note that the shape of \eqref{E.cubic} is such that the sum of the three roots in $z$ of the resulting equation with $u=0$ vanishes. Therefore if two intersections are constant, the third one also is.

\medskip

Summing up, with these normalizations we have achieved that for $s\in S'_c$, the transformation $L_s^{-1}$ fixes   the three points of $E^*$ on $U=0$ and  fixes the line $U=0$; this entails that it fixes pointwise the line $U=0$. 
Also, it  fixes the flex $(0:1:1)$.

Then, up to a nonzero factor,  $L_s$, for $s\in S'_c$  has rows of the  shape  $(1,0,0), (0,a,0), (0,a-1,1)$, for some $a=a_s$. 
However now we use that $L_s$ sends the shape \eqref{E.cubic} into the same shape, i.e. with vanishing coefficient of $z^2$. This yields $a=1$, so $L_s$ is the identity and all $\tilde X_s$ for $s\in S'_c$, would be isomorphic. This is impossible, as we are assuming that only finitely many $\tilde X_s$ are isomorphic to a given one.

This proves the assertion.
\end{proof}

NOTE: Note that the proof yields also that  our  section cannot be a (torsion or non-torsion) constant on a subscheme with iso-constant fibers, if the base has dimension $2$. (It makes sense to speak of a `constant' section if the scheme is isotrivial.)

\medskip


\begin{remark} The above argument allows us to be explicit as to the family with iso-constant elliptic curves. Note that after base change over $S_c$ we obtain a trivial square-elliptic family $(E^*)^2\times R_c$, where $R_c$ is  a finite cover of $S_c$,  and where  the given section corresponds to a map from $R_c$ to $(E^*)^2$. 

This map will  now be shown to be not dominant, and actually with image a certain curve that we now describe. And the fibers will all be curves in $R_c$. The  argument now shows that we can obtain these curves as follows.

 First, by a translation on $u=U/T$ we can assume that the section is obtained by taking two differences  (in the elliptic-group sense) between the three points of the curve on the line $U=0$. 

Now take a given one elliptic  curve of the family, denoted $E^*$, having a flex at $(0:1:1)$, taken as origin,  and passing through $(1:0:1)$, where the section (to $(E^*)^2$) can be assumed to correspond to two of the points on the line $U=0$. (Note that there will be three such points if the point at infinity is non-special, and the sum of the three will be zero.) As in the argument for Proposition \ref{P.nt}, this curve $E^*$ can be gotten through a projective transformation fixing the isomorphism class of the curve $\tilde X$ of genus $2$ leading to the said elliptic curve (as in the description above). 

 Now take the projective transformations $L_s$ fixing both the flex  $(0:1:1)$ and the point $(1:0:1)$. 
 Their matrices have rows of the shape  $(a,b,-b), (c,d,-c), (e+b-a, e+c-b, e)$. They form a group $\Gamma$ of (projective) dimension $4$.   This groups acts faithfully on the projective plane $\P_2$, and also on the set of equations (for instance if it fixes an equation it induces an automorphism of an elliptic curve, with two fixed points, thus must induce the identity). 
 
 \medskip
 

 Now there is a subvariety $V\subset \Gamma$ 
 formed by the $L\in \Gamma$ such that the equation for $L(E^*)$ has zero coefficients of $Z^2$ (there are two coefficients). Note that 
 $V$ has  dimension $\ge 2$, as a nonempty  fiber (above the point $(0,0)\in\A^2$) of a map $\Gamma\to\A^2$.
 
  But then the dimension is exactly $2$. Indeed it cannot have dimension $\ge 3$,  because otherwise there would exist inside $V$ an irreducible  curve corresponding to a single isomorphism class of $\tilde X$. 
Then the corresponding polynomials $P^2-Q^3$ would all be related by a homography, up to a constant factor; but then the same would be true individually of the pairs $P,Q$ (because there are only finitely many pairs producing a given $P^2-Q^3$, as will appear in the mentioned forthcoming paper. \footnote{This may be proved directly as follows. 
Let $F(u)=P^2(u)-Q^3(u)$, where $F$ is fixed and  without multiple roots.  If there are infinitely many solutions in polynomials $P,Q$ of degrees resp. $3,2$, then there is an algebraic family of solutions of dimension $1$. In other words we may assume that the coefficients of $P,Q$ are rational functions on a certain curve $Z$, not all constant. Ler $\delta$ be a derivation on $\C(Z)$ with constant field $\C$, and apply it to the equation, obtaining $2(\delta P)P=3(\delta Q)Q^2$. But $P,Q$ must be coprime, hence $Q^2$ divides $\delta P$ and $P$ divides $\delta Q$, forcing $\delta P,\delta Q$ to vanish. So after all they are constant.})

 In conclusion, we would have a $1$-parameter family of matrices in $\Gamma$ carrying a fixed  cubic $z^3-3Q(u)z-2P(u)$ into a constant multiple of $z^3-3\rho^2(cu+d)^2Q(\gamma(u))z-2\rho^3(cu+d)^3P(\gamma(u))$ for a nonzero $\rho$ and a $\gamma(u)=(au+b)/(cu+d)\in PGL_2$. This would produce a $1$-dimensional family of  automorphisms $(u,z)\mapsto (\gamma(u),z/\rho(cu+d))$ of the corresponding elliptic curve, which is impossible (e.g. by general theory, or note that these transformation would have to stabilise the set of flexes, or else the branch points of the map $u:E^*\to\P_1$). 
 
 Naturally it is possible to prove the contention about $\dim V$ also by explicit calculation. (Or else we could also use the above Proposition; in fact, if $V$ has dimension $3$, we would get a subvariety of dimension $2$ such that the section is constant on it, contrary to the last remark in the proof of the Proposition.)


 
 \medskip
 

Then,  the  cubic curves of the family  may be taken just the $E_L:=L(E^*)$ for $L\in V$: they are all isomorphic, but they come from generically non-isomorphic curves $\tilde X$.

This also shows that the image of the section, from the above space $R_c$, and represented in $E^*\times E^*$, is of the shape $x_c\times E^*$ for a certain $x_c\in E^*$ (depending also on the family of points at infinity, but now we are assuming this is given).  This $x_c$ is represented by the point $(0:1:0)$ in $E^*$ (with the chosen flex as origin). 

When $x_c$ is not torsion, the method does not yield anything on the space $S_c$. On the contrary, when it is torsion, we have pencils like in the example presented in \S\ref{S.pencil}. 

This description also yields a rather complete picture of the whole. 

In the paper in progress (also with Lombardo) we plan to provide full detail for this issue too, and to present all types of examples that can be deduced from it. 
\end{remark}






\medskip


\bigskip

\address{
Pietro Corvaja\\
Dipartimento di Matematica e Informatica\\
Universit\`a di Udine\\
Via delle Scienze, 206\\
33100 Udine  \\
Italy
}
{pietro.corvaja@dimi.uniud.it}

\address{
Umberto Zannier\\
Scuola Normale Superiore\\
Piazza dei Cavalieri, 7\\
56126 Pisa \\
Italy
}
{umberto.zannier@sns.it}


\bigskip





 

\begin{thebibliography}{99}


\smallskip


\bibitem{ACZ} - Y. Andr\'e, P. Corvaja, U. Zannier (with an appendix by Z. Gao), The Betti map for a section of an abelian scheme, Inv. Math. 2020.

\smallskip

\bibitem{BerGyo} - A. Berczes, J-H. Evertse, K. Gy\"ory, Effective results for Diophantine equations over finitely generated domains, Acta Arithmetica 163 (2014), 71-100.

\smallskip

\bibitem{Bi} - Yu. Bilu, Effective analysis of integral points on algebraic curves,
Israel J. Math. 90 (1995), 235--252.

\smallskip

\bibitem{BiI} - Yu. Bilu, M. Illengo, Effective Siegel's theorem for modular curves,
Bull. London Math. Soc. 43 (2011), 673--688.


\smallskip


\bibitem{BG} - E. Bombieri, W. Gubler, Heights in Diophantine Geometry, Cambridge Univ, Press, 2006.

\smallskip

\bibitem{CMZ} - P. Corvaja, D. Masser, U. Zannier, Torsion hypersurfaces on abelian scemes and Betti coordinates, Math. Annalen, 2017.


 
\smallskip

\bibitem{DZ} - R. Dvornicich, U. Zannier, Fields containing values of algebraic functions, Annali Sc. Normale Sup. Cl. Sci,  21 (1994), 421-443.

\smallskip

\bibitem{GGW} - B. Grechuk, T. Grechuk, A. Wilcox, Diophantine equations with three monomials, preprint 2023.

\smallskip

\bibitem{LP} -  A. Landesman, B. Poonen, Obbstructions to applying the Baker-Bilu method for determining integral points on curves, preprint 2023.  


\smallskip

\bibitem{Lomb} D. Lombardo, E. L. Garcia, Ch. Ritzenthaler, J. Sijsling, Decomposing Jacobians Via Galois covers, Experimental Math. 32 (2023), 218-240.

\smallskip




\smallskip

\bibitem{S} - A. Schinzel, Polynomials with special regard to reducibility, Cambridge Univ., Press, 1999.

\smallskip

\bibitem{Se} - J-P. Serre, Lectures on the Mordell-Weil Theorem, Aspects of Mathematics, Vieweg, third edition 2013.

\smallskip


\smallskip

\bibitem{Z2} - U. Zannier, An effective solution of a certain diophantine problem, Rend. Sem. Mat. Univ. Padova, 93 (1995), 177--183.

\smallskip

\bibitem{Z} - U. Zannier, Some problems of Unlikely Intersections in Arithmetic and Geometry, Annals of Math. Studies 181, Princeton U.P., 2012.

\end{thebibliography}
\end{document}